\RequirePackage{ifpdf}
\ifpdf % We are running pdfTeX in pdf mode
\documentclass[pdftex]{sigma}
\else
\documentclass{sigma}
\fi

\begin{document}

\allowdisplaybreaks

\renewcommand{\PaperNumber}{054}

\FirstPageHeading

\ShortArticleName{Horizontal Forms of Chern Type on Complex Finsler Bundles}

\ArticleName{Horizontal Forms of Chern Type\\ on Complex Finsler Bundles}

\Author{Cristian IDA}

\AuthorNameForHeading{C.~Ida}

\Address{Department of Algebra, Geometry and Dif\/ferential Equations,\\
Transilvania University of Bra\c{s}ov, Str. Iuliu Maniu 50, Bra\c{s}ov 500091, Rom\^{a}nia}
\Email{\href{mailto:cristian.ida@unitbv.ro}{cristian.ida@unitbv.ro}}

\ArticleDates{Received October 28, 2009, in f\/inal form June 30, 2010;  Published online July 09, 2010}

\Abstract{The aim of this paper is to construct horizontal Chern forms of a holomorphic vector bundle using complex Finsler structures. Also, some properties of these forms are studied.}


\Classification{53B40; 57R20}

\section{Introduction and preliminaries}

In \cite{Fa} J.J.~Faran posed an open question: \textit{Is it possible to define Chern forms of a holomorphic vector bundle using complex Finsler structures?} Recently, in \cite{Ai2, C-W, Ch-W} and \cite{Wo} was studied the Chern form of the associated hyperplane line bundle in terms of complex Finsler structures and some of its applications.

The main purpose of this note is to obtain  horizontal Chern forms of a complex Finsler bundle, following the general construction of Chern forms \cite{G-H, Ko1, Va}. Firstly, we def\/ine a $h''$-cohomology with respect to the horizontal conjugated dif\/ferential operator and using the partial Bott complex connection~\cite{Ai1}, we construct horizontal invariant forms of Chern type on the total space. Also, we prove that the corresponding horizontal classes of Chern type do not depend on the complex Finsler structure chosen in a family given by conformal changes of horizontal type. Next, it is proved that these forms live on the associated projectivized bundle and are invariant by a linear family of complex Finsler structures given by inf\/initesimal deformations.

Let $\pi:E\rightarrow M$ be a holomorphic vector bundle of rank $m$ over a $n$-dimensional complex manifold $M$. Let $(U,(z^{k}))$, $k=1,\dots ,n$ be a local chart on $M$ and $s=\{s_{a}\}$, $a=1,\dots ,m$ a~local frame for the sections of $E$ over $U$. It is well known that this chart canonically induces another one on $E$ of the form $(\pi^{-1}(U),u=(z^{k},\eta^{a}))$, $k=1,\dots ,n$, $a=1,\dots ,m$, where $s=\eta^{a}s_{a}$ is a~section on $E_{z}=\pi^{-1}(z)$, for all $z\in M$. If $(\pi^{-1}(U'),u'=(z^{\prime k},\eta^{\prime a}))$ is another local chart on~$E$, then the transition laws of these coordinates are
\[
z^{\prime k}=z^{\prime k}(z),\qquad \eta^{\prime a}=M^{a}_{b}(z)\eta^{b},
\]
where $M^{a}_{b}(z)$, $a,b=1,\dots, m$ are holomorphic functions on $z$ and $\det M^{a}_{b}\neq0$.

As we already know \cite{Ai1, Ai2}, the total space $E$ has a natural structure of $m+n$ dimensional complex manifold because the transition functions $M^{a}_{b}$ are holomorphic. We consider the complexif\/ied tangent bundle $T_{\mathbb{C}}E$ of the real tangent bundle $T_{\mathbb{R}}E$ and its decomposition $T_{\mathbb{C}}E=T'E\oplus T''E$, where $T'E$ and $T''E=\overline{T'E}$ are the holomorphic and antiholomorphic tangent bundles of $E$, respectively. The vertical holomorphic tangent bundle $V'E=\ker \pi_{*}$ is the relative tangent bundle of the holomorphic projection $\pi$. A local frame f\/ield on $V'_{u}E$ is $\{\partial/\partial\eta^{a}\}$, $a=1,\dots ,m$ and $V'_{u}E$ is isomorphic to the sections module of~$E$ in~$u$.

A supplementary subbundle of $V'E$ in $T'E$, i.e.\ $T'E=H'E\oplus V'E$ is called a \textit{complex nonlinear connection} on $E$, brief\/ly c.n.c. A local basis for the horizontal distribution $H'_{u}E$, called \textit{adapted} for the c.n.c., is $\{\delta/\delta z^{k}=\partial/\partial z^{k}-N^{a}_{k}\partial/\partial\eta^{a}\}$, $k=1,\dots ,n$, where $N^{a}_{k}(z,\eta)$ are the local coef\/f\/icients of the c.n.c. In the sequel we consider the abreviate notations: $\partial_{k}=\partial/\partial z^{k}$, $\delta_{k}=\delta/\delta z^{k}$ and $\dot{\partial}_{a}=\partial/\partial\eta^{a}$. The adapted basis denoted by $\{\delta_{\overline{k}}=\delta/\delta\overline{z}^{k}\}$ and $\{\dot{\partial}_{\overline{a}}=\partial/\partial\overline{\eta}^{a}\}$, for $H''E=\overline{H'E}$ and $V''E=\overline{V'E}$ distributions are obtained respectively by conjugation anywhere. Also, we note that the adapted cobases are given by $\{dz^{k}\}$, $\{\delta\eta^{a}=d\eta^{a}+N^{a}_{k}dz^{k}\}$, $\{d\overline{z}^{k}\}$ and $\{\delta\overline{\eta}^{a}=d\overline{\eta}^{a}+N^{\overline{a}}_{\overline{k}}d\overline{z}^{k}\}$ which span the dual bundles $H^{\prime *}E$, $V^{\prime *}E$, $H^{\prime\prime *}E$ and $V^{\prime\prime *}E$, respectively.

\begin{definition}[\cite{Ko2}] A strictly pseudoconvex complex Finsler structure on $E$ is a positive real valued function $F:E\rightarrow \mathbb{R}_{+}\cup\{0\}$ which satisf\/ies the following conditions:
\begin{enumerate}\itemsep=0pt
\item[1)] $F^{2}$ is smooth on $E^{*}=E\setminus \{0\}$;

\item[2)] $F(z,\eta)\geq 0$ and $F(z,\eta)=0$ if and only if $\eta=0$;

\item[3)] $F(z,\lambda\eta)=|\lambda|F(z,\eta)$ for any $\lambda\in\mathbb{C}^{*}=\mathbb{C}\setminus \{0\}$;

\item[4)] the complex hessian $(h_{a\overline{b}})=(\dot{\partial}_{a}\dot{\partial}_{\overline{b}}(F^{2}))$ is positive def\/inite and determines a~hermitian metric tensor on the f\/ibers of the vertical bundle $V'E \setminus  \{{\rm zero\ section}\}$.
\end{enumerate}
\end{definition}
\begin{definition}
The pair $(E,F)$ is called a complex Finsler bundle.
\end{definition}

According to \cite{A-P, Ai2, Mu1}, a c.n.c.\ on $(E,F)$ depending only on the complex Finsler structure $F$ is the Chern--Finsler c.n.c., locally given by
\begin{gather}
\stackrel{\rm CF}{N^{a}_{k}}=h^{\overline{c}a}\partial_{k}\dot{\partial}_{\overline{c}}\big(F^{2}\big),
\label{I2}
\end{gather}
where $(h^{\overline{c}a})$ is the inverse of $(h_{a\overline{c}})$.

An important property of the Chern--Finsler c.n.c.\ (see \cite{Ai2, Mu1}), is
\begin{gather}
[\delta_{j},\delta_{k}]=\big(\delta_{k}\stackrel{\rm CF}{N_{j}^{a}}-\delta_{j}\stackrel{\rm CF}{N_{k}^{a}}\big)\dot{\partial}_{a}=0.
\label{I3}
\end{gather}

Throughout this paper we consider adapted frames and coframes with respect to the Chern--Finsler c.n.c.

According to \cite{Ai1}, a partial connection $\nabla$ of $(1,0)$ type on $V'E$, def\/ined by
\begin{gather*}
\nabla_{X}Y=v'[X,Y]
%\label{I4}
\end{gather*}
for any $X\in\Gamma(H'E)$ and $Y\in\Gamma(V'E)$,
is called \textit{the partial Bott complex connection} of $(E,F)$.
Here $v'$ is the natural vertical projector.

The $(1,0)$-connection form $\omega=(\omega^{a}_{b})$ of $\nabla$ is locally given by
\begin{gather}
\omega^{a}_{b}=L^{a}_{bk}dz^{k}, \qquad L^{a}_{bk}=\dot{\partial}_{b}\big(\stackrel{\rm CF}{N^{a}_{k}}\big)
\label{I5}
\end{gather}
and the $(1,1)$-horizontal curvature form $\Omega=(\Omega^{a}_{b})$ of $\nabla$ is locally given by
\begin{gather}
\Omega^{a}_{b}=\Omega^{a}_{bj\overline{k}}dz^{j}\wedge d\overline{z}^{k}, \qquad \Omega^{a}_{bj\overline{k}}=-\delta_{\overline{k}}(L^{a}_{bj}).
\label{I6}
\end{gather}

\begin{remark}
We notice that the partial Bott complex connection is in fact the classical Chern--Rund connection, f\/irst introduced by H.~Rund~\cite{Ru} in the case when  $E=T'M$ is the holomorphic tangent bundle of a complex manifold $M$.
\end{remark}

\subsection*{$\boldsymbol{h''}$-cohomology}

Similar to \cite{P-M}, let us consider  the set $\mathcal{A}^{p,q,r,s}(E)$ of all $(p,q,r,s)$-forms with complex values on~$E$, locally def\/ined by
\[
\varphi=\frac{1}{p!q!r!s!}\varphi_{I_{p}\overline{J_{q}}A_{r}\overline{B_{s}}}dz^{I_{p}}\wedge d\overline{z}^{J_{q}}\wedge\delta\eta^{A_{r}}\wedge\delta\overline{\eta}^{B_{s}},
\]
where $I_{p}=(i_{1}\dots i_{p})$, $J_{q}=(j_{1}\dots j_{q})$, $A_{r}=(a_{1}\dots a_{r})$, $B_{s}=(b_{1}\dots b_{s})$ and these forms can be nonzero only when they act on $p$ vectors from $\Gamma(H'E)$, on $q$ vectors from $\Gamma(H''E)$, on $r$ vectors from $\Gamma(V'E)$ and on $s$ vectors from $\Gamma(V''E)$, respectively.

By (\ref{I3}), we get the following decomposition of the exterior dif\/ferential
\begin{gather*}
d\mathcal{A}^{p,q,r,s}(E)\subset\mathcal{A}^{p+1,q,r,s}(E)
\oplus\mathcal{A}^{p,q+1,r,s}(E)\oplus\mathcal{A}^{p,q,r+1,s}(E)\oplus\mathcal{A}^{p,q,r,s+1}(E)\oplus{}
\\
\quad{}\oplus\mathcal{A}^{p+1,q+1,r-1,s}(E)\oplus\mathcal{A}^{p+1,q,r-1,s+1}(E)\oplus\mathcal{A}^{p+1,q+1,r,s-1}(E)\oplus\mathcal{A}^{p,q+1,r+1,s-1}(E)
\end{gather*}
which allows us to def\/ine eight morphisms of complex vector spaces if we consider dif\/ferent components in the above decomposition.

If we consider  the subset $\mathcal{A}^{p,q,0,0}(E)\subset\mathcal{A}^{p,q,r,s}(E)$ of all horizontal forms of $(p,q,0,0)$ type, which can be nonzero only when they act on $p$ vectors from $\Gamma(H'E)$ and $q$ vectors from $\Gamma(H''E)$, then
\[
d\mathcal{A}^{p,q,0,0}(E)\subset\mathcal{A}^{p+1,q,0,0}(E)\oplus\mathcal{A}^{p,q+1,0,0}(E)\oplus\mathcal{A}^{p,q,1,0}(E)\oplus\mathcal{A}^{p,q,0,1}(E).
\]

In particular, similar to \cite{Zho1, Zho2}, the horizontal dif\/ferential operators are def\/ined by
\[
d^{\prime h}: \ \mathcal{A}^{p,q,0,0}(E)\rightarrow\mathcal{A}^{p+1,q,0,0}(E), \qquad d^{\prime \prime h}: \ \mathcal{A}^{p,q,0,0}(E)\rightarrow\mathcal{A}^{p,q+1,0,0}(E)
\]
where, for any $\varphi=\frac{1}{p!q!}\varphi_{I_{p}\overline{J_{q}}}dz^{I_{p}}\wedge d\overline{z}^{J_{q}}\in\mathcal{A}^{p,q,0,0}(E)$, we have
\[
(d^{\prime h}\varphi)_{i_{1}\dots i_{p+1}\overline{J_{q}}}=\sum_{k=1}^{p+1}(-1)^{k-1}\delta_{i_{k}}\big(\varphi_{i_{1}\dots \widehat{i_{k}}\dots i_{p+1}\overline{J_{q}}}\big)
\]
and
\begin{gather}
(d^{\prime\prime h}\varphi)_{I_{p}\overline{j_{1}}\dots \overline{j_{q+1}}}=(-1)^{p}\sum_{k=1}^{q+1}(-1)^{k-1}\delta_{\overline{j_{k}}}\big(\varphi_{I_{p}\overline{j_{1}}\dots \widehat{\overline{j_{k}}}\dots \overline{j_{q+1}}}\big).
\label{I9}
\end{gather}

We note that by (\ref{I3}), we have $d^{\prime\prime h}\circ d^{\prime\prime h}=0$. Thus, with respect to the operator $d^{\prime\prime h}$ we can def\/ine  \textit{the $h''$-cohomology groups} of $(E,F)$ for $(p,q,0,0)$-forms by
\begin{gather}
H^{q}\big(E,\Phi^{p,0,0}\big)=\frac{{\rm ker} \{d^{\prime\prime h}:\mathcal{A}^{p,q,0,0}\rightarrow\mathcal{A}^{p,q+1,0,0}\}}{{\rm Im}\, \{d^{\prime\prime h}:\mathcal{A}^{p,q-1,0,0}\rightarrow\mathcal{A}^{p,q,0,0}\}},
\label{I10}
\end{gather}
where $\Phi^{p,0,0}$ is the sheaf of germs of $(p,0,0,0)$-forms $d^{\prime\prime h}$-closed.

\section{Horizontal forms of Chern type}

In this section, using the horizontal curvature form of partial Bott complex connection, we construct horizontal forms of Chern type on a complex Finsler bundle.

With the notations from the previous section, we have $\omega^{a}_{b}\in\mathcal{A}^{1,0,0,0}(E)$ and $\Omega^{a}_{b}\in\mathcal{A}^{1,1,0,0}(E)$. Also, by ({\ref{I6}}) and ({\ref{I9}}) we have $\Omega=d^{\prime\prime h}\omega$ and taking into account the relation $d^{\prime\prime h}\circ d^{\prime\prime h}=0$, it follows that
$\Omega$ is a $d^{\prime\prime h}$-closed dif\/ferential form.

Let $gl(m,\mathbb{C})\approx\mathbb{C}^{m \times m}$ be the Lie algebra of the linear general group $Gl(m,\mathbb{C})$. A \textit{symmetric polynomial} $f\in S^{j}(Gl(m,\mathbb{C}))$ is invariant if
\[
f\big(M^{-1}X_{1}M,\dots ,M^{-1}X_{j}M\big)=f(X_{1},\dots ,X_{j})
\]
for every $M\in Gl(m,\mathbb{C})$ and $X_{j}\in gl(m,\mathbb{C})$.

Also, it is well known \cite{G-H, Ko1, Va}, that the algebra of invariant symmetric polynomials on $gl(m,\mathbb{C})$ is generated by the elementary symmetric functions $f_{j}\in S^{j}(Gl(m,\mathbb{C}))$ given by
\begin{gather}
\det\left(I_{m}-\frac{1}{2\pi i}X\right)=\sum_{j=0}^{m}f_{j}(X)=1-\frac{1}{2\pi i}\,{\rm tr}\,X+\cdots +\left(-\frac{1}{2\pi i}\right)^{m}\det X.
\label{II2}
\end{gather}

Next, we give a known result, but we express it in a form adapted to our study.
\begin{proposition}
At local changes on $E$, the horizontal curvature of partial Bott complex connection changes by the rule
\[
\Omega'=M^{-1}\Omega M,
\]
where $M=(M^{a}_{b}(z))$ are the transition functions of the holomorphic bundle $E$.
\end{proposition}
\begin{proof}
By \cite{Mu1}, we know that the change rules of the adapted frames and coframes are
\begin{gather}
\delta_{j}=\frac{\partial z^{\prime k}}{\partial z^{j}}\delta'_{k},\qquad \dot{\partial}_{b}=M^{a}_{b}\dot{\partial}'_{a}, \qquad dz^{j}=\frac{\partial z^{j}}{\partial z^{\prime k}}dz^{\prime k}, \qquad \delta\eta^{b}=M^{b}_{a}\delta\eta^{\prime a}
\label{II4}
\end{gather}
and its conjugates.

On the other hand (see 7.1.9 in \cite{Mu1}), the local coef\/f\/icients of the Chern--Finsler c.n.c.\ change by the rule
\begin{gather}
\frac{\partial z^{\prime k}}{\partial z^{j}}\stackrel{\rm CF}{N^{\prime a}_{k}}=M^{a}_{b}\stackrel{\rm CF}{N^{b}_{j}}-\frac{\partial M^{a}_{b}}{\partial z^{j}}\eta^{b}.
\label{II5}
\end{gather}
Finally, it follows by (\ref{I5}), (\ref{I6}), (\ref{II4}) and (\ref{II5}) that $\Omega^{\prime d}_{c}=M^{b}_{c}\Omega^{a}_{b}M^{d}_{a}$.
\end{proof}

Now, we consider the $(j,j,0,0)$-invariant forms $\mathcal{C}^{h}_{j}(\nabla)$ def\/ined by
\[
\det \left(I_{m}-\frac{1}{2\pi i}\Omega\right)=\sum_{j=0}^{m}\mathcal{C}^{h}_{j}(\nabla).
\]

Then, it follows by (\ref{II2})
\[
\mathcal{C}^{h}_{j}(\nabla)=\frac{(-1)^{j}}{(2\pi i)^{j}j!}\sum_{a,b}\delta^{b_{1}\dots b_{j}}_{a_{1}\dots a_{j}}\Omega^{a_{1}}_{b_{1}}\wedge\dots \wedge\Omega^{a_{j}}_{b_{j}},
\]
where $\delta^{b_{1}\dots b_{j}}_{a_{1}\dots a_{j}}$ are the Kronecker symbols.

Because $d^{\prime\prime h}$ satisf\/ies the property
\[
d^{\prime\prime h}(\varphi\wedge\psi)=d^{\prime\prime h}\varphi\wedge\psi+(-1)^{\deg\varphi}\varphi\wedge d^{\prime\prime h}\psi
\]
for any $\varphi\in\mathcal{A}^{p,q,0,0}(E)$ and $\psi\in\mathcal{A}^{p',q',0,0}(E)$, then by $d^{\prime\prime h}\Omega=0$ we obtain that the dif\/ferential forms $\mathcal{C}^{h}_{j}(\nabla)$ are $d^{\prime\prime h}$-closed.
\begin{definition}
$\mathcal{C}^{h}_{j}(\nabla)$ are called horizontal forms of Chern type of order $j$ of the complex Finsler bundle $(E,F)$.
\end{definition}

By (\ref{I10}) these forms def\/ine the horizontal cohomology classes
\[
c^{h}_{j}(\nabla)=[\mathcal{C}^{h}_{j}(\nabla)]\in H^{j}\big(E,\Phi^{j,0,0}\big)
\]
which are called \textit{horizontal classes of Chern type of order $j$} of $(E,F)$. In particular, the f\/irst horizontal class of Chern type is represented by the $(1,1,0,0)$-form
\[
\mathcal{C}^{h}_{1}(\nabla)=-\frac{1}{2\pi i}\Omega^{a}_{aj\overline{k}}dz^{j}\wedge d\overline{z}^{k}.
\]

\begin{proposition}
The first horizontal class of Chern type $c^{h}_{1}(\nabla)$ is invariant by conformal change  of horizontal type: $F^{2}\mapsto e^{\sigma(z)}F^{2}$.
\end{proposition}
\begin{proof}
Let $F^{2}\mapsto\widetilde{F}^{2}:=e^{\sigma(z)}F^{2}$ be a conformal change of complex Finsler structures on $E$, where $\sigma(z)$ is a smooth function on $M$. Because $\widetilde{h}_{a\overline{b}}=e^{\sigma(z)}h_{a\overline{b}}$, using (\ref{I2}) we get that the Chern--Finsler c.n.c. changes by the rule
\[
\stackrel{\rm CF}{\widetilde{N}^{a}_{k}}=\stackrel{\rm CF}{N^{a}_{k}}+\frac{\partial\sigma}{\partial z^{k}}\eta^{a}.
\]

Thus, the $(1,0)$-connection form of $\nabla$ changes by the rule
\begin{gather}
\widetilde{\omega}=\omega+d^{\prime h}\sigma\otimes I,
\label{II13}
\end{gather}
where $I$ is the identity endomorphism of $V'E$.

Applying $d^{\prime \prime h}$ in the (\ref{II13}) relation, we get the change rule of the $(1,1)$-horizontal curvature form of $\nabla$, namely
\[
\widetilde{\Omega}=\Omega+d^{\prime\prime h}d^{\prime h}\sigma\otimes I
\]
which ends the proof.
\end{proof}

In fact the above proposition remains valid for all higher horizontal Chern classes, namely:
\begin{proposition}
The horizontal classes of Chern type $c_{j}^{h}(\nabla)$ are invariant by conformal change  of horizontal type: $F^{2}\mapsto e^{\sigma(z)}F^{2}$ for any $j=1,\dots ,m$.
\end{proposition}
\begin{proof}
It is classical that for a square matrix $A$ of order $p$ we have
\[
\det(A-\lambda I)=\sum_{j=0}^{p}(-1)^{j}\Delta_{p-j}(A)\lambda^{j},
\]
where
\[
\Delta_{j}(A)=\frac{1}{j!}\sum_{\alpha,\beta}\delta^{\beta_{1}\dots \beta_{j}}_{\alpha_{1}\dots \alpha_{j}}a^{\alpha_{1}}_{\beta_{1}} \cdots  a^{\alpha_{j}}_{\beta_{j}}
\]
is the sum of the principal minors of order $j$ of the matrix $A$ (see for instance \cite[p.~235]{Va}).

Now, taking into account $d^{\prime\prime h}\Delta_{j}(\Omega)=0$, by direct calculations we have
\begin{gather*}
\mathcal{C}^{h}_{j}(\widetilde{\nabla})  =  \frac{(-1)^{j}}{(2\pi i)^{j}}\Delta_{j}(\Omega+d^{\prime\prime h}d^{\prime h}\sigma\otimes I)\\
\phantom{\mathcal{C}^{h}_{j}(\widetilde{\nabla})}  =  \frac{(-1)^{j}}{(2\pi i)^{j}}\left(\Delta_{j}(\Omega)+\frac{1}{j!}\sum_{k=1}^{j}\Delta_{j-k}(\Omega)\wedge(d^{\prime \prime h}d^{\prime h}\sigma)^{k}\right)\\
\phantom{\mathcal{C}^{h}_{j}(\widetilde{\nabla})} = \mathcal{C}^{h}_{j}(\nabla)+d^{\prime\prime h}\left(\frac{(-1)^{j}}{(2\pi i)^{j}j!}\sum_{k=1}^{j}\Delta_{j-k}(\Omega)\wedge d^{\prime h}\sigma\wedge(d^{\prime \prime h}d^{\prime h}\sigma)^{k-1}\right)
\end{gather*}
which says that $\mathcal{C}^{h}_{j}(\widetilde{\nabla})$ and $\mathcal{C}^{h}_{j}(\nabla)$ are in the same $d^{\prime \prime h}$-cohomology class.
\end{proof}

It is well known that for every $\varphi\in\mathcal{A}^{p,q}(M)$, $\varphi^{*}=\varphi\circ\pi\in\mathcal{A}^{p,q,0,0}(E)$ is the natural lift of $\varphi$ to the total space $E$.
\begin{remark}
If the complex Finsler structure $F$ comes from a hermitian structure $h$ on $E$, namely $F^{2}(z,\eta)=h(\eta,\eta)=h_{a\overline{b}}(z)\eta^{a}\overline{\eta}^{b}$, then the coef\/f\/icients of the partial Bott complex connection~$\nabla$ are independent of the f\/iber coordinates $\eta^{a}$, and these coincides with the classical connection coef\/f\/icients in hermitian bundles. In this case it is easy to see that $\mathcal{C}_{j}^{h}(\nabla)$ coincide with the classical Chern forms of order $j$ of hermitian bundle $(E,h)$ lifted to the total space $E$.
\end{remark}

In the sequel we suppose that the base manifold $M$ is compact. Note that there is a natural $\mathbb{C}^{*}=\mathbb{C}\setminus \{0\}$ action on $E^{*}=E\setminus \{0\}$ and the associated projectivized bundle is def\/ined by $PE=E^{*}/\mathbb{C}^{*}$ with the projection $p:PE\rightarrow M$. Each f\/iber $P_{z}(E)=P(E_{z})$ is isomorphic to the $(m-1)$-dimensional complex projective space $P^{m-1}(\mathbb{C})$. The pull-back bundle $p^{-1}E$ is a~holomorphic vector bundle of rank $m$ over $PE$. Thus, the local complex coordinates~$(z,\eta)$ on~$E$ may be also considered as a local complex coordinates system for $PE$ as long as $\eta^{1},\dots ,\eta^{m}$ is considered as a homogeneous coordinate system for f\/ibers. All the geometric objects on $E$ which are invariant after replacing $\eta$ by $\lambda\eta$, $\lambda\in\mathbb{C}^{*}$ are valid on~$PE$. The reason for working on~$PE$ rather than on~$E$ is that $PE$ is compact if~$M$ is compact (see~\cite{Ko2, Ko3}).

Now, we simply denote by $\mathcal{A}^{p,q}(PE)$ the set of all horizontal forms of $(p,q,0,0)$ type on $E$ whose coef\/f\/icients are zero homogeneous with respect to f\/iber coordinates, namely $\varphi_{I_{p}\overline{J_{q}}}(z,\lambda\eta)=\varphi_{I_{p}\overline{J_{q}}}(z,\eta)$, for any $\lambda\in\mathbb{C}^{*}$.

By the homogeneity conditions of the complex Finsler structure $F$ (see \cite{Ai1, Ai2, Mu1}), the local coef\/f\/icients of the Chern--Finsler c.n.c. are given by
\begin{gather}
\stackrel{\rm CF}{N^{a}_{k}}=L^{a}_{bk}\eta^{b}
\label{II15}
\end{gather}
and, moreover
\begin{gather}
\stackrel{\rm CF}{N^{a}_{k}}(z,\lambda\eta)=\lambda\stackrel{\rm CF}{N^{a}_{k}}(z,\eta),\qquad \forall\,\lambda\in\mathbb{C}^{*}.
\label{II16}
\end{gather}

Then, it follows by (\ref{II15}) and (\ref{II16})
\[
L^{a}_{bk}(z,\lambda\eta)=L^{a}_{bk}(z,\eta),\qquad \forall\,\lambda\in\mathbb{C}^{*}
\]
which says that the $(1,0)$-connection forms $\omega^{a}_{b}$ of $\nabla$ live on $PE$. Similarly, by (\ref{I6}) we get $\Omega^{a}_{bj\overline{k}}(z,\lambda\eta)=\Omega^{a}_{bj\overline{k}}(z,\eta)$ for any $\lambda\in\mathbb{C}^{*}$ which says that the $(1,1)$-curvature forms $\Omega^{a}_{b}$ of $\nabla$ live on $PE$. Thus,  $\mathcal{C}^{h}_{j}(\nabla)\in\mathcal{A}^{j,j}(PE)$.

Finally, we note that properties of horizontal Chern forms in relation with basic properties of classical Chern forms must be studied, as well as the independence of these forms with respect to some families of complex Finsler structures on a holomorphic vector bundle. Another important problem to solve is to describe the obstructions corresponding to these classes. Here we are able to respond only partially to these questions.

Applying some results from \cite{Ai}, we get the invariance of horizontal Chern forms by a family of complex Finsler structures given by inf\/initesimal deformations.

We consider a $1$-parameter family $\{F^{2}_{t}\}$, $t\in\mathbb{R}$ of pseudoconvex complex Finsler structures on a holomorphic vector bundle $E$. Let us put $F^{2}_{0}=F^{2}$ and let $\nabla_{0}=\nabla$ its partial Bott complex connection. The inf\/initesimal deformation $V$ induced by $F^{2}_{t}$ is def\/ined by $V:=\left(\frac{\partial F^{2}_{t}}{\partial t}\right)|_{t=0}$ and its components with respect to a f\/ixed frame f\/ield $s=\{s_{a}\}$, $a=1,\dots ,m$ of $E$ are given by $V_{a\overline{b}}=\left(\frac{\partial h_{t\,a\overline{b}}}{\partial t}\right)|_{t=0}$. We put $V^{a}_{b}=h^{\overline{c}a}V_{b\overline{c}}$ and we consider it as a section of the bundle ${\rm End}(V'E)$. If we consider $F^{2}_{t}=F^{2}+tV$ for suf\/f\/iciently small $t$ so that $F^{2}_{t}$ remains pseudoconvex then in Theorem~5.2 from~\cite{Ai} it is proved that if $\nabla V^{a}_{b}=0$ then $\nabla_{t}=\nabla$ and $\Omega_{t}=\Omega$. Here $\nabla Z=(d^{\prime h}Z^{a}+Z^{b}\omega^{a}_{b})\dot{\partial}_{a}$, for any $Z=Z^{a}\dot{\partial}_{a}\in\Gamma(V'E)$. Thus, we can conclude
\begin{proposition}
If $\nabla V^{a}_{b}=0$ then the horizontal Chern forms $\mathcal{C}_{j}^{h}(\nabla)$ of $(E,F)$ are invariant by the linear family of complex Finsler structures given by $F^{2}_{t}=F^{2}+tV$.
\end{proposition}

We notice that in Aikou's paper \cite{Ai}, the partial connection is considered in the pull-back bundle $p^{-1}E$ over $PE$ and the calculations are similar as on $V'E$.

\subsection*{Acknowledgements}

The author is grateful to the anonymous
referees and would like to thank them for generous suggestions and comments. Also, I warmly thank Professor Gheorghe Piti\c{s} for fruitful conversations concerning this topics.

\pdfbookmark[1]{References}{ref}

\LastPageEnding

\end{document}